\documentclass{gtart}
\def\ifplaintex{\expandafter\ifx\csname documentclass\endcsname\relax}

\ifplaintex 
\else
\fi

\expandafter\ifx\csname beginpicture\endcsname\relax
\expandafter\ifx\csname documentclass\endcsname\relax
\input pictex \else
\input prepictex \input pictex \input postpictex \fi\fi

\def\gt{{\mathsurround=0pt\it $\cal G\mskip-2mu$eometry \&\ 
$\cal T\!\!$opology}}        %  journal title in recommended style

\def\gtp{{\mathsurround=0pt\it $\cal G\mskip-2mu$eometry \&\ 
$\cal T\!\!$opology $\cal P\!$ublications}}  % GT publications

\def\lognumber#1{\def\thelognumber{#1}}
\def\volumenumber#1{\def\thevolumenumber{#1}}
\def\papernumber#1{\def\thepapernumber{#1}}
\def\volumeyear#1{\def\thevolumeyear{#1}}

\def\pagenumbers#1#2{\def\startpage{#1}\def\finishpage{#2}}
\def\published#1{\def\publishdate{#1}}
\def\proposed#1{\def\theproposer{#1}}
\def\seconded#1{\def\theseconders{#1}}
\def\received#1{\def\receiveddate{#1}}
\def\revised#1{\def\reviseddate{#1}}
\def\accepted#1{\def\accepteddate{#1}}

\def\asciiaddress#1{\def\theasciiaddress{#1}}

\long\def\asciiabstract#1{\long\def\theasciiabstract{#1}}

\let\\\par\let\thelognumber\relax
\let\thevolumenumber\relax\let\thepapernumber\relax
\let\thevolumeyear\relax\let\thesamplenumber\relax\let\startpage\relax
\let\finishpage\relax\let\publishdate\relax\let\receiveddate\relax
\let\reviseddate\relax\let\accepteddate\relax\let\theasciititle\relax
\let\theasciiauthors\relax\let\theasciiaddress\relax
\let\theasciiabstract\relax
\let\theasciiemail\relax\let\theshortauthors\relax\let\theshorttitle\relax

\long\def\maketitlep{   % start of definition of \maketitlep

\count0=\startpage

\gt\hfill      %   Journal title (top left) 
\beginpicture
\setcoordinatesystem units <0.33truein, 0.33truein> point at 2.2 0.9
\setplotsymbol ({$\cal G$})
\plotsymbolspacing=9truept
\circulararc 315 degrees from 0 1 center at 0 0
\setplotsymbol ({$\cal T$})
\circulararc 315 degrees from 1 -1 center at 1 0
\endpicture
\break
{\small\ifx\thesamplenumber\relax % sample?  
Volume \else Sample
\fi\thevolumenumber\ (\thevolumeyear)
\startpage--\finishpage\nl
Published: \publishdate}
\vglue 0.5truein plus 0.4fil minus 0.1truein

{\parskip=0pt\leftskip 0pt plus 1fil\def\\{\par\smallskip}{\ifplaintex\large
\else\Large\fi\bf\thetitle}\par\medskip}   

\vglue 0pt plus 0.1fil 

{\parskip=0pt\leftskip 0pt plus 1fil\def\\{\par}{\sc\theauthors}
\par\medskip}

\vglue 0pt plus 0.1fil 

{\small\parskip=0pt\let\newline\\
{\leftskip 0pt plus 1fil\def\\{\par}{\sl\theaddress}\par}
\expandafter\ifx\theemail\relax    % email address?
\relax\else\vglue 5pt plus 0.02fil minus 2pt\def\\{\stdspace{\rm 
and}\stdspace} 
\cl{Email:\stdspace\tt\theemail}\fi
\ifx\theurl\relax                  % URL given?
\relax\else\vglue 5pt plus 0.02fil minus 2pt\def\\{\stdspace{\rm 
and}\stdspace}
\cl{URL:\stdspace\tt\theurl}\fi\par}

\vglue 7pt plus 0.3fil minus 3pt

{\bf Abstract}
\vglue 5pt plus 0.1fil minus 2pt

\theabstract

\vglue 7pt plus 0.3fil minus 3pt

{\bf AMS Classification numbers}\quad Primary:\quad \theprimaryclass

Secondary:\quad \thesecondaryclass

\vglue 5pt plus 0.3fil minus 2pt

{\bf Keywords}\quad \thekeywords

\vglue 10pt plus 0.5fil minus 5pt

{\small  Proposed: \theproposer\hfill Received: \receiveddate\nl
Seconded: \theseconders\hfill 
\ifx\reviseddate\relax                         % paper revised?
Accepted: \accepteddate                        % no
\else
Revised: \reviseddate                          % yes
\fi}
\eject
}       %  end of definition of \maketitlep

\let\maketitlepage\maketitlep
\let\maketitle\maketitlepage

\font\phead=cmsl9 scaled 950
\font=cmsl9 scaled 1050
\font\pnum=cmbx10 scaled 913
\font\lnum=cmbx10 
\font\pfoot=cmsl9 scaled 950
\font=cmsl9 scaled 1050
\ifplaintex
\headline{\vbox to 0pt{\vskip -4.5mm\line{\small\phead\ifnum
\count0=\startpage ISSN 1364-0380 (on line)
1465-3060 (printed) \hfill {\pnum\folio}\else\ifodd\count0\def\\{ }% 
\ifx\theshorttitle\relax\thetitle\else\theshorttitle\fi\hfill{\pnum\folio}
\else\def\\{ and }{\pnum\folio}\hfill\ifx\theshortauthors\relax\theauthors
\else\theshortauthors\fi\fi\fi}\vss}}
\footline{\vbox to 0pt{\vglue 0mm\line{\small\pfoot\ifnum\count0=\startpage
\copyright\ \gtp\hfill\else
\gt, Volume \thevolumenumber\ (\thevolumeyear)\hfill\fi}\vss
}}
\else
\makeatletter
\def\@oddhead{{\small\lhead\ifnum\count0=\startpage ISSN 1364-0380 (on line)
1465-3060 (printed) \hfill {\lnum\number\count0}\else\ifodd\count0
\def\\{ }\ifx\theshorttitle\relax \thetitle \else\theshorttitle\fi\hfill
{\lnum\number\count0}\else\def\\{ and }{\lnum\number\count0}
\hfill\ifx\theshortauthors\relax 
\theauthors\else\theshortauthors\fi\fi\fi}}\def\@evenhead{\@oddhead}
\def\@oddfoot{\small\lfoot\ifnum\count0=\startpage\copyright\ \gtp\hfill\else
\gt, Volume \thevolumenumber\ (\thevolumeyear)\hfill\fi}
\def\@evenfoot{\@oddfoot}
\makeatother
\fi

\newwrite\gtoutfile
\long\gdef\makeheadfile{  %%% start of definition of \makeheadfile
{\def\\{, }\def\s{ }
\immediate\openout\gtoutfile head.xxx
\immediate\write\gtoutfile{Proxy-for: \ifx\theasciiauthors\relax
\theauthors\else\theasciiauthors\fi\s<\ifx\theasciiemail\relax\theemail\else\theasciiemail\fi>}
\immediate\write\gtoutfile{\noexpand\\}
\immediate\write\gtoutfile{Authors: \ifx\theasciiauthors\relax
\theauthors\else\theasciiauthors\fi}
{\def\\{ }\immediate\write\gtoutfile{Title: \ifx\theasciititle\relax
\thetitle\else\theasciititle\fi}}
\immediate\write\gtoutfile{Subj-class: GT or SG or MG etc}
\immediate\write\gtoutfile{MSC-class: \theprimaryclass\ifx\thesecondaryclass\relax\else, \thesecondaryclass\fi}
\immediate\write\gtoutfile{Journal-ref: Geom. Topol. \thevolumenumber
(\thevolumeyear) \startpage-\finishpage}
\immediate\write\gtoutfile{Comments: Published by Geometry and Topology at}
\immediate\write\gtoutfile{\s\s http://www.maths.warwick.ac.uk/gt/GTVol\thevolumenumber/paper\thepapernumber.abs.html}
\immediate\write\gtoutfile{\noexpand\\}
\immediate\write\gtoutfile{}
\ifx\theasciiabstract\relax
\immediate\write\gtoutfile{\theabstract}\else
\immediate\write\gtoutfile{\theasciiabstract}\fi
\immediate\write\gtoutfile{}
\immediate\write\gtoutfile{\noexpand\\}
\immediate\write\gtoutfile{}
\immediate\closeout\gtoutfile}}  %%% end of definition of \makeheadfile

\def\maketitlepage{\maketitlep\makeheadfile}
\let\maketitle\maketitlepage

\lognumber{310}
\volumenumber{7}\papernumber{21}\volumeyear{2003}
\pagenumbers{757}{771}
\received{28 March 2003}
\revised{26 October 2003}
\published{30 October 2003}
\accepted{29 October 2003}
\proposed{Benson Farb}
\seconded{Leonid Polterovich, Joan Birman}

\usepackage{amsmath,amssymb}

\newtheorem{thm}{Theorem}[section]
 \newtheorem{lemma}[thm]{Lemma}
\newtheorem{cor}[thm]{Corollary}

\newtheorem{prop}[thm]{Proposition}

\theoremstyle{definition}
\newtheorem{defn}[thm]{Definition}

\newtheorem{ex}[thm]{Example}

\newtheorem{remark}[thm]{Remark}

\DeclareMathOperator{\Fix}{Fix}
\DeclareMathOperator{\Diff}{Diff}

\DeclareMathOperator{\Per}{Per}
\newcommand{\R}{\mathbb R}

\newcommand{\Z}{\mathbb Z}

\def\ti{\tilde}
\def\sinfty{S_{\infty}}
\def\slthreez{SL(3, \mathbb Z)}
\def\ti{\tilde}
\def\sinfty{S_{\infty}}

\def\G{{\cal G}}

\def\A{{\mathbb A}}

\begin{document}
\title{Area preserving group actions on surfaces}
\author{John Franks\\Michael Handel}
\address{Department of Mathematics, Northwestern 
University\\Evanston, IL 60208-2730, USA}

\secondaddress{Department of Mathematics, CUNY, Lehman 
College\\Bronx, NY 10468, USA}

\asciiaddress{Department of Mathematics, Northwestern 
University\\Evanston, IL 60208-2730, USA\\and\\Department 
of Mathematics, CUNY, Lehman 
College\\Bronx, NY 10468, USA}

\email{john@math.northwestern.edu, handel@g230.lehman.cuny.edu}
\begin{abstract}
Suppose $\G$ is an almost simple group containing a subgroup isomorphic
to the three-dimensional integer Heisenberg group.  For example any finite
index subgroup of $\slthreez$ is such a group.  The main result of this
paper is that every action of $\G$ on a closed oriented surface by
area preserving diffeomorphisms factors through a finite group.
\end{abstract}
\asciiabstract{Suppose G is an almost simple group containing a
subgroup isomorphic to the three-dimensional integer Heisenberg group.
For example any finite index subgroup of SL(3,Z) is such a group.
The main result of this paper is that every action of G on a closed
oriented surface by area preserving diffeomorphisms factors through a
finite group.}


\maketitlepage
\section{Introduction and notation}

This article is motivated by the program of classifying actions of
higher rank lattices in simple Lie groups on closed manifolds.
More specifically,  we are concerned here with   actions of  $SL(n, \Z)$ with $n \ge 3$,  
on closed oriented surfaces.
A standard example of such an action is given by
projectivizing the usual action of $\slthreez$ on $\R^3.$ Our objective
here is to show that there are essentially no such actions which
are symplectic or, what amounts to the same thing in this case, area
preserving.  

R. Zimmer conjectured (see Conjecture 2 of \cite{Z2}) that any
$C^\infty$ volume preserving action of a finite index subgroup of
$SL(n,\Z)$ on a closed manifold with dimension less than $n$, factors
through an action of a finite group.  We prove this in the case that
the dimension of the manifold is $2.$  L Polterovich \cite{P} previously
provided a proof of this result for surfaces of genus at least one.
Our aim was to extend this to the case of the sphere.  As it happens our
techniques, which are quite different from his, are equally applicable to
any genus so we present the argument in full generality.

\begin{defn}\label{defn: almost simple}
A group $\G$ is called {\it almost simple} if every normal subgroup is
either finite or has finite index.
\end{defn}

The Margulis normal subgroups theorem (see Theorem IX.5.4 of
\cite{Mar} or 8.1.2 of \cite{Z1}) asserts that an irreducible lattice
in a semi-simple Lie group with $\R-$rank $\ge 2$ is almost simple.
In particular, any finite index subgroup of $SL(n, \Z)$ with $n \ge 3$
is almost simple.

Suppose $S$ is a closed oriented surface and $\omega$ is a smooth
volume form.  We will generally assume a fixed choice of $\omega$ and
refer to a diffeomorphism $F\co  S \to S$ which preserves $\omega$ as an
{\em area preserving diffeomorphism}.  We denote the group of
diffeomorphisms preserving $\omega$ by $\Diff_\omega(S)$ and its
identity component by $\Diff_\omega(S)_0.$ Equivalently
$\Diff_\omega(S)_0$ is the group of diffeomorphisms which preserve
$\omega$ and are isotopic to $id$.

A key ingredient (and perhaps limitation) of our approach to this
problem is the fact that a finite index subgroup of $SL(n, \Z)$ with
$n \ge 3$ always contains a subgroup isomorphic to the
three-dimensional integer Heisenberg group.  Recall that this is the
group of upper triangular $3 \times 3$ integer matrices with all
diagonal entries equal to $1$.  Our main result is:

\begin{thm}  \label{no actions} 
Suppose $\G$ is an almost simple group containing a subgroup
isomorphic to the three-dimensional integer Heisenberg group, e.g.\ any finite index subgroup of $SL(n, \Z)$ with $n \ge 3$.  If $S$ is a
closed oriented surface then every homomorphism $\phi\co  \G \to
\Diff_\omega(S)$ factors through a finite group.
\end{thm}

Given $\G, S$ and $\phi$ as in Theorem~\ref{no actions}, let $\phi' \co 
\G \to MCG(S)$ be the homomorphism to the mapping class group $MCG(S)$
induced by $\phi$. Since nilpotent subgroups of $MCG(S)$ are virtually
abelian (see \cite{BLM}), the kernel $\G_0$ of $\phi'$ has an infinite order
element in the integer Heisenberg subgroup.  Hence by almost simplicity $\G_0$
has finite index in $\G$.  Moreover, $\G_0$ contains a finite index
subgroup of the integer Heisenberg subgroup and hence contains a subgroup
isomorphic to the three-dimensional integer Heisenberg group. Thus in proving
Theorem~\ref{no actions} there is no loss in replacing $\G$ with
$\G_0$.  In other words, we may assume that the image of $\phi$ is
contained in $\Diff_\omega(S)_0$.

If $n \ge 4$, then any analytic action by $SL(n, \Z)$ on a closed
oriented surface $S$ factors through a finite group. This was proved
by Farb and Shalen \cite{FS} for $S \ne T^2$ and by Rebelo \cite{Re}
for $S = T^2$.  (Farb and Shalen proved this for $S=T^2$ under the
assumption that the action is area preserving.)  Polterovich (see
Corollary 1.1.D of \cite{P}) proved that if $n \ge 3$, then any action
by $SL(n, \Z)$ on a closed surface $S$ other than $S^2$ and $T^2$ by
area preserving diffeomorphisms factors through a finite group.  All
of these results can be stated in greater generality than we give
here.  There is also an analogous result of D.~Witte (\cite{Wi}),
which asserts that a homomorphism $\phi\co  SL(n, \Z) \to Homeo(S^1)$
must factor through a finite group if $n \ge 3.$

We are grateful to Benson Farb for introducing us to the problem and for several very helpful converations.  We are also grateful to David Fisher for suggestions that significantly streamlined and helped organize the paper.

\section{Hyperbolic structures and normal form}

Some of our proofs rely on mapping class group techniques that use
hyperbolic geometry.  In this section we establish notation and recall a result from \cite{fh:periodic}.  For further details see, for example, \cite{casble:nielsen}).

   Let $S$ be
a closed orientable surface.  We will say that a connected open subset
$M$ of $S$ has {\em negative Euler characteristic} if $H_1(M, \R)$ is
infinite dimensional or if $M$ is of finite type and the usual
definition of Euler characteristic has a negative value.  If $M$ has negative Euler characteristic then $M$ supports a complete hyperbolic structure.  

We use the Poincar\'e disk model for the hyperbolic plane $H$.  In
this model, $H$ is identified with the interior of the unit disk and
geodesics are segments of Euclidean circles and straight lines that
meet the boundary in right angles. A choice of complete hyperbolic structure on
$M$ provides an identification of the universal cover $\ti M$ of $M$
with $H$.  Under this identification covering translations become
isometries of $H$ and geodesics in $M$ lift to geodesics in $H$.  The
compactification of the interior of the unit disk by the unit circle
induces a compactification of $H$ by the \lq circle at infinity\rq\
$\sinfty$.  Geodesics in $H$ have unique endpoints on $\sinfty$.
Conversely, any pair of distinct points on $\sinfty$ are the endpoints
of a unique geodesic.

 Suppose that
$F\co  S \to S$ is an orientation preserving homeomorphism of a closed surface $S$ and  that $M
\subset S$ is an open connected $F$-invariant  set with negative Euler characteristic.  Equip $M$ with a complete hyperbolic structure and let  $f =F|_M \co M\to M$. We use the identification of $H$ with $\ti M$ and write $\ti f \co  H \to
H$ for lifts of $f\co  M \to M$ to the universal cover.  A fundamental
result of Nielsen theory is that every lift $\ti f \co  H \to H$ extends
uniquely to a homeomorphism (also called) $\ti f \co  H \cup \sinfty \to
H \cup \sinfty$.  (A proof of this fact appears in Proposition 3.1 of
\cite{han:fpt}).  If $f \co  M \to M$ is isotopic to the identity then
there is a unique lift $\ti f$, called the {\it identity lift,} that
commutes with all covering translations and whose extension over
$\sinfty$ is the identity.

Every covering translation $T \co  H \to H$ extends to a homeomorphism
(also called) $T \co  H \cup \sinfty \to H \cup \sinfty$.  If $\gamma
\subset M$ is a closed geodesic and $\ti \gamma \subset H$ is a lift
to $H$ then there is an extended covering translation $T$ whose only
fixed points $T^+$ and $T^-$ are the endpoints of $\ti \gamma$. If
$f(\gamma)$ is isotopic to $\gamma$ then there is a lift $\ti f \co  H
\cup \sinfty \to H \cup \sinfty$ that fixes $T^{\pm}$ and commutes
with $T$.  The quotient space of $H \cup (\sinfty \setminus T^{\pm})$
by the action of $T$ is a closed annulus on which $\ti f$ induces a
homeomorphism denoted $\hat f \co  \A \to \A$.

The following result is an immediate corollary of Theorem 1.2 and Lemma 6.3 of \cite{fh:periodic}. If $\Fix(F)$ is finite, then it is just a special case of the Thurston classification theorem \cite{Th}.  

\begin{thm} \label{can form} 
Suppose that $F\co  S \to S$ is a diffeomorphism of an orientable
closed surface, that $F$ is isotopic to the identity and that $\Fix(F) \ne \emptyset$.  Then  
  there is a finite set $R$ of simple closed curves in $ S
\setminus \Fix(F)$ and a homeomorphism $\phi$ isotopic to $F$ rel $\Fix(F)$
such that:
\begin{description}
\item [(1)] $\phi$ setwise fixes disjoint open annulus neighborhoods $\A_j
\subset (S \setminus \Fix(F))$ of the elements $\gamma_j \in R $.  The
restriction of $\phi$ to $cl(\A_i)$ is a non-trivial Dehn twist of a closed annulus.
\end{description}

Let $\{S_i\}$ be the components of $S \setminus \cup \A_j$, let $X_i =
\Fix(F) \cap S_i$, and if $X_i$ is finite, let $M_i = S_i \setminus X_i$.

\begin{description}
\item [(2)] If $X_i$ is infinite then $\phi|_{S_i}$= identity. 
\item [(3)] If $X_i$ is finite, then   $M_i$ has negative Euler characteristic and $\phi|_{M_i}$  is either pseudo-Anosov  or the identity.  
\end{description}
\end{thm}

We say that $\phi$ is a {\it normal form for $F$} and that $R$ is {\it
the set of reducing curves for $\phi$}.
If  $R$ has minimal cardinality among all sets of reducing curves for all normal forms for $F$, then we say that {\it $R$ is a minimal set of reducing curves}.

\section{The proof of Theorem~\ref {no actions}}

We denote by $[G,H]$ the commutator
$G^{-1}H^{-1}GH.$

\begin{prop} \label{commutator} 
If $G,H \co  S \to S$ are diffeomorphisms that are isotopic to the identity and that
commute with their commutator $F=[G,H]$ then     $F$ 
  is isotopic to the identity rel $\Fix(F)$.
\end{prop}

\begin{remark} \label{iterate} In any group containing elements
$G$ and $H$, if $F = [G,H]$ commutes with $G$ and $H$ then 
it is easy to see that $[G^k, H] = [G,H]^k = [G,H^k]$ for any
$k \in \Z.$  Hence $F^{k_1k_2} = [G^{k_1},H^{k_2}]$ for all $k_1, k_2$.
\end{remark}

As an immediate consequence of this and Proposition~\ref{commutator}
we have 

\begin{cor} \label{cor: commutator power} 
Suppose that $G,H \co S \to S$ are diffeomorphisms that are isotopic to the identity and   that commute with their commutator
$F = [G,H]$. Then for all $n>0$, $F^n$ is isotopic to the identity rel $\Fix(F^n)$.

\end{cor}

Before proving Proposition~\ref{commutator} we state and prove some required lemmas.
Denote the closed annulus by $\A$ and its boundary components by $\partial_0 \A$ and $\partial_1 \A$.  The universal cover $\ti \A$ is identified with $\mathbb R \times [0,1]$.

\begin{lemma} \label{annulus comm}  
Assume that $u_i \co  \A \to \A$, $i=0,1$, are homeomorphisms that
preserve $\partial_0 \A$ and $\partial_1 \A$, that $u_i$ commutes with
$v=[u_1,u_2]$ and that $\ti u_i \co \ti \A \to \ti \A$ are lifts of
$u_i$.  Then $\ti v = [\ti u_1, \ti u_2]$ fixes at least one point in
both $\partial_0 \ti \A$ and $\partial_1 \ti \A$.
\end{lemma}

\begin{proof} 
Let $p_1 \co \ti \A \to {\mathbb R}$ be projection onto the first
coordinate and let $T \co  \ti \A \to \ti \A$ be the indivisible covering
translation $T(x,y) = (x+1,y)$. We claim that $|p_1(\ti v(\ti x)) -
p_1(\ti x)| < 8$ for all $\ti x \in \partial \ti \A$.

Since the lift $\ti u_i$ commutes with $T$, $\ti v$ is independent of the
choice of $\ti u_i$.  We may therefore assume that $|p_1(\ti u_i(\ti
x_i)) -p_1(\ti x_i)| < 1 $ for some $\ti x_i \in \partial_0 \ti \A$
and hence that $|p_1(\ti u_i(\ti x)) -p_1(\ti x)| < 2$ for all $\ti x
\in \partial_0 \ti \A$.  The claim for $\ti x \in \partial_0 \ti \A$
follows immediately.  The analogous argument on $\partial_1 \ti \A$
completes the proof of the claim.

The preceding argument holds for any $u_i$ that preserve the
components of $\partial \A$ and in particular for all iterates
$u_i^N$.  We conclude (Remark~\ref{iterate}) that $|p_1(\ti
v^{N^2}(\ti x)) - p_1(\ti x)| < 8$ for all $\ti x \in \partial \ti \A$
and all $N$.  Since the restriction of $\ti v$ to each boundary
component of $\ti \A$ is an orientation preserving homeomorphism of
the line, these homeomorphisms both have points with a bounded forward
orbit, and any such homeomorphisms of the line must have a fixed point.
\end{proof}

An isolated end $E$ of an open set $U \subset S$ has
neighborhoods of the form $N(E) = S^1 \times [0,1)$.  The set $fr(E) =
cl_S(N(E)) \setminus N(E)$, called the {\it frontier of $E$}, is 
independent of the choice of $N(E)$.  We will say that an end is
{\em trivial} if $fr(E)$ is a single point.

\begin{lemma}  \label{lem: extension}  
If $E$ is a non-trivial isolated end of an open subset $U$ of $S$ then
there is a manifold compactification of $E$ by a circle $C$ satisfying
the following property.  If $F$ is any orientation preserving
homeomorphism of $U \cup fr(E)$ and if $U$ is $F$-invariant then
$F|_U$ extends to a homeomorphism $f$ of $U \cup C$. Moreover, if $G$
is another orientation preserving homeomorphism of $U \cup fr(E)$ and
if $g$ is the extension of $G|_U$ over $U \cup C$ then:
\begin{description}
\item [1)]if $F$ is isotopic to $G$ relative to $fr(E)$ then $f$ is isotopic to $g$ relative to  $C$.  In particular, if $F$ is isotopic to the identity relative to $fr(E)$ then $f$ is isotopic to the identity relative to $C$.
\item [2)] $fg$ is the extension of $(FG)|_U$.
\end{description}

If $E$ is a trivial end with $fr(E) = \{x\}$ and 
$F$ and $G$ are local diffeomorphisms on a neighborhood of $x$, 
there is a manifold compactification of $E$ by a circle $C$ and  extensions $f,g$ to homeomorphisms of $U \cup C$ satisfying property 2).

\end{lemma} 

\begin{proof}  Assume at first that $E$ is non-trivial.
The existence of $C$ and $f$ is a consequence of the theory of prime
ends (see \cite{M2} for a good modern exposition).  Since $C$ is a
boundary component of $U \cup C$, the extension $f$ is the unique
extension of $F|_U$ over $C$.  Property (2) follows immediately.

    If $F$ and $G$ are as (1), then $FG^{-1}$
is isotopic to the identity relative to $fr(E)$.  It suffices to show that $fg^{-1}$ is isotopic to the identity relative to $C$ since precomposing with $g$ then gives the desired isotopy of $f$ to $g$.  We may therefore assume that $G$ is the identity.

Given  an isotopy $H_t$ of $F$ to the identity relative to $fr(E)$, extend $H_t|_U$ by the identity on $C$ to define   $h_t \co  U \cup C \to U \cup C$. It suffices to show that $h_t(x)$ varies continuously in both $t$ and $x$. This is clear if $x \in U$ so suppose that $x \in C$.  Choose a disk neighborhood system $\{W_i\}$  for $x$.   It suffices to show that for all $i$ and $t$ there exists $j$ so that $h_s(W_j \cap U) \subset W_i \cap U$ for all $s$ sufficiently close to $t$.  The frontier of $W_i$ is an arc $\sigma_i$ that intersects $C$ exactly in its endpoints.  It suffices to show that $h_s(\sigma_j)$ is disjoint from $\sigma_i$ and lies in the same component of $U \setminus \sigma_i$ as does $\sigma_j$.   

The system $\{W_i\}$ can be chosen with three important properties (see \cite{M2} for details.) First, the interior of $\sigma_i$, thought of as an open arc in $U$, is the interior of a closed arc $\sigma_i'$ in  $U \cup fr(E)$ with endpoints in $fr(E)$. Second,  the paths $\sigma_i'$s are mutually disjoint and converge to a point  $z \in fr(E)$. The  third property is that for all $i$ and $t$, the interior of $H_t(\sigma_i')$, thought of as an open arc in $U$ is the interior of a closed arc in $U \cup C$ with the same endpoints as $\sigma_i$; this closed arc is by definition $h_t(\sigma_i)$ .   

     Fix $i$ and $t$.  Since $H_t(z) = z$,  there exists $j > i$ such that  $H_t(\sigma_j') \cap \sigma_i' = \emptyset$. By compactness of the closed arcs, $H_s(\sigma_j') \cap \sigma_i' = \emptyset$ for all $s$ sufficiently close to $t$.  Since $h_s(\sigma_j)$ and $\sigma_j$ have the same endpoints, $h_s(\sigma_j)$ lies in the same component of $U \setminus \sigma_i$ as does $\sigma_j$. This completes the proof in the case that $E$ is non-trivial.

In the case that $E$ is trivial and $F,G$ are local diffeomorphisms
we can construct $C$ by blowing up the point $x$ to obtain $C$ and
the continuous extensions to $U \cup C.$  The blowing up construction
is functorial so property 2) is satisfied.
\end{proof}

\noindent{\bf Proof of Proposition~\ref{commutator}}\qua   We may assume without loss that $\Fix(F) \ne \emptyset$.  Since $G$ and $H$ commute with $F$, $\Fix(F)$ is $G$-invariant and $H$-invariant. Let $\phi$ be a canonical form for $F$ with minimal reducing set $R$ and let $X_i, S_i$ and $M_i$ be as in   Theorem~\ref{can form}.  It suffices to show that if $X_i$ is finite then $\phi|_{M_i}$ is not pseudo-Anosov and that   $R = \emptyset$.  These properties are unchanged if $F$ is replaced by an iterate so there is no loss in replacing $G$,  $H$ and $F$  by $G^n$, $H^n$ and $F^{n^2}$ for some $n > 0$.  Lemma 6.2 of \cite{fh:periodic} implies that $G$ and $H$ permute the elements of $R$ up to isotopy relative to Fix($F$).  We may therefore assume that $G$ and $H$ fix the elements of $R$ up to isotopy relative to Fix($F$).

  We first rule out the possibility that some $\phi |_{M_i}$ is
pseudo-Anosov.    
For each $M_i \subset M$ there are
well defined elements $<F|_{M_i}>,<G|_{M_i}>$ and
$<H|_{M_i}>$ in the mapping class group of $M_i$ defined by \lq
restricting\rq\ $F, G$ and $H$ to $M_i$.  For concreteness
we give the argument for $G$.  Since $G$ preserves the isotopy
class of each element of $R$, $G$ preserves the isotopy class of
each component of $\partial M_i$ and $G \simeq G_1$ where $G_1(M_i)
= M_i$.  If $G_2$ is another such homeomorphism then, since $G_1
\simeq G_2$, $G_1|_{M_i} \simeq G_2|_{M_i}$ are isotopic as
homeomorphisms of $M_i$. Define $<G|_{M_i}>$ to be the isotopy class
determined by $G_1|_{M_i}$.

The isotopy classes $<F|_{M_i}>,<G|_{M_i}>$ and $<H|_{M_i}>$
are determined by the actions of $F, G$ and $H$ on isotopy
classes of simple closed curves in $M_i$. Thus $<F|_{M_i}> =
[<G|_{M_i}><H|_{M_i}>]$ commutes with both $<G|_{M_i}>$ and
$<H|_{M_i}>$ for all $i$.

If some $<F|_{M_i}> =<\phi|_{M_i}>$ is pseudo-Anosov with expanding
lamination $\Lambda$, then $<G|_{M_i}>$ and $<H|_{M_i}>$ are
contained in the stabilizer of $\Lambda$.  But this stabilizer is
virtually Abelian (see, for example, Lemma 2.3 of
\cite{han:commuting}) in contradiction to the fact that
$<F|_{M_i}>^{m^2} = [<G|_{M_i}>^m <H|_{M_i}>^m]$ is non-trivial for
all $m$.  This completes the proof that $\phi |_{M_i}$ is not
pseudo-Anosov. Thus $\phi$ is the identity on the complement of the annuli $\A_j$.

If $R \ne \emptyset$ choose $\gamma \in R$, let $U$ be  
  the component of $S \setminus \Fix(F)$ that contains $\gamma$ and let $f, g$ and $h$ be $F|_U,G|_U$ and $H|_U$ respectively.   There are three cases to consider.  The first is that $U$ is an open annulus.   By Lemma~\ref{lem: extension} there is a compactification of $U$ to a closed annulus $\A$ and homeomorphisms $\hat f,\hat g, \hat h, \hat \phi \co  \A \to \A$ that respectively extend $ f,g,h, \phi|U \co  U \to U$ and that satisfy 
\begin{description}
\item [(1)]  $\hat f$ commutes with both $\hat g$ and $\hat h$. 
\item [(2)]  $\hat f = [\hat g, \hat h]$. 
\item  [(3)]  $\hat f$ is isotopic rel $\partial \A$ to $\hat \phi$.
\end{description}
  By hypothesis, 
\begin{description}
\item[(4)] $\hat \phi$ is isotopic rel $\partial \A$ to a non-trivial Dehn twist.
\end{description}
so 
\begin {description}
\item [(5)]  if $\hat \alpha$ is an arc with endpoints on distinct components of $\partial \A$, then $\hat f(\hat \alpha)$ is not isotopic rel endpoints to $\alpha$.
\end{description}
  Property (5) contradicts Lemma~\ref{annulus comm} and so completes the proof in this first case.

We may now assume that $U$ has negative Euler characteristic and hence supports a complete hyperbolic structure.  The second case is that $\gamma$ is not peripheral in $U$.  There is no loss in assuming that $\gamma$ is a geodesic.    Choose a lift $\ti \gamma \subset H$ to the universal cover of $U$ and let
$T\co  H \cup \sinfty \to H \cup \sinfty$ be the extended indivisible covering
translation that preserves $\ti \gamma$.  Choose lifts $ \ti g, \ti h
\co  H \cup \sinfty \to H \cup \sinfty$ of $g,h$ that commute with $T$ and so fix
the endpoints $T^\pm$ of $\ti \gamma$ in $\sinfty$.  Then $\ti f
=[\ti g, \ti h]$ is a lift of $f$ that commutes with $T$ and there is a lift $\ti \phi$ of $\phi|_U$ that is equivariantly isotopic   to $\ti f$.     Let $\A$ be the annulus obtained as the quotient
space of $H \cup (\sinfty \setminus T^{\pm})$ by the action of $T$, and let $\hat \phi,\hat f, \hat
g, \hat h \co  \A \to \A$ be the homeomorphisms induced by $\ti \phi, \ti f, \ti g$
and $\ti h$.  Items (2) and (3) above are immediate. Since $\ti f \ti g$ and $\ti g \ti f$ project to the
same homeomorphism of $U$ and commute with $T$, they differ by an
iterate of $T$.  Thus $\hat f\hat g =\hat g \hat f$.  The symmetric
argument shows that $\hat f$ commutes with $\hat h$ so (1) is satisfied.  In this  case the fixed point set of $\hat \phi$ intersects each  component of $\partial \A$ in a Cantor set so we must replace (4) with a weaker property (see the proof of Proposition 7.1 of \cite{fh:periodic}) for details):
\begin{description}
\item[($4'$)] $\hat \phi$ is isotopic to a non-trivial Dehn twist relative to a closed set  that intersects both components of $\partial \A$.
\end{description}

Property (5) follows from ($4'$) so the  proof concludes as in the previous case.  

   The last case is that $\gamma$ is peripheral in $U$ and is a minor variation on the second case.   By Lemma~\ref{lem: extension}, the end of $U$ corresponding to $\gamma$ can be compactified by a circle and the homeomorphisms $ f,g,h, \phi|_U$ can be extended to homeomorphisms of the resulting space $U_*$. There is a hyperbolic structure on $U_*$ in which $\gamma$ is isotopic to a peripheral geodesic $\gamma_*$.  The universal cover $\ti U_*$ is naturally identified with the interior of the convex hull in $H$ of a Cantor set $C \subset \sinfty$ and  is compactifed by a circle consisting of the union of $C$ with the full pre-image of $\gamma_*$. (See for example page 175 of \cite{han:nielsen}.) The proof now proceeds exactly as in the second case using $\ti U_*$ and its circle compactification in place of $H$ and $\sinfty$.     
\qed

\begin{lemma} \label{lem: hyperbolic complements} 
Suppose $G,H \co  S \to S$ are area preserving, orientation preserving
diffeomorphisms that are isotopic to the identity and that commute with their commutator $F=[G,H].$ Then
except in one case all components of $S \setminus \Fix(F)$ have negative Euler characteristic. The one exception is the case that $S = S^2$
and $\Fix(F)$ consists of exactly two points.  In this case all
components of $S \setminus \Fix(F^2)$ have negative Euler characteristic. 
\end{lemma}

\begin{proof}
Let $U$ be a component of $S \setminus \Fix(F)$.  By \cite{brn-kis} it
is $F$-invariant.  The Poincar\'e recurrence theorem and the Brouwer plane translation theorem imply that every area preserving homeomorphism of the open disk has a fixed point; thus $U$ can not be  an open disk.  If $S = T^2$, then the mean rotation vector of $F$ is zero since $F$ is a commutator; by \cite{CZ}, $F$ has fixed points and   $U \ne T^2$.  
We are left only with the case that $U$ is an open annulus.
By Lemma~\ref{lem: extension} we can compactify $U$ to a closed annulus
$\A$ and extend $F$ to $\A$ continuously.  Suppose first
that we are not in the exceptional case that $S = S^2$ and $\Fix(F)$ consists
of exactly two points.  Then by Lemma~\ref{lem: extension} the map $F$ is 
the identity on at least one component of $\partial \A.$

Because $G$ and $H$ are area preserving and preserve $\Fix(F)$ there
exist $k,l >0$ such that $G^k(U) = U$ and $H^l(U) = U$.  By Lemma~\ref{lem: extension} we can extend $G^k,H^l$
 to $\A$ and by doubling $k$ and $l$ if necessary we may
assume $G^k,H^l$ preserve the boundary components of $\A$.

Then $F^{kl} = [G^k, H^l]$ by Remark~\ref{iterate}.  Let $\ti \A$ be the
universal covering space of $\A$ and let $u, v\co  \ti \A \to \ti \A$ be
lifts of $G^k$ and $H^l$ respectively.  Then $w = [u,v]$ is a lift
of $F^{kl}$.  Let $\ti F$ be a lift of $F$ which is the identity on
one boundary component of $\ti \A.$   Then $\ti F^{kl}$ is also a
lift of $F^{kl}$.  By Lemma~\ref{annulus comm} the map $w$ has
fixed points in both boundary components of $\ti \A.$  It follows
that $w$ and $\ti F^{kl}$ have a common fixed point and hence they
are equal since they are both lifts of the same map.

But the mean rotation of the lift $w$ is zero since it is a
commutator.  It follows from Theorem 2.1 of \cite{Fr6}, that it has an
interior fixed point.  This implies that $\ti F$ has an interior
periodic point and hence an interior fixed point by the Brouwer plane
translation theorem.  This, in turn, implies $F$ has a fixed point in
$U$ which is a contradiction.

We are left with the single exceptional case that $S = S^2$ and
$\Fix(F) = \{p,q\}$  so $U$ is the open annulus $S^2 \setminus \{p,q\}$.
By Lemma~\ref{lem: extension} we can compactify
$U$ to form an annulus $\A$ and  extend $F,G$ and $H$ to
orientation preserving, area preserving homeomorphisms of $\A$.  
Then $G^2$ and $H^2$  must preserve the boundary components of $\A$.
Suppose first that in addition one of $G$ and $H$ (say $G$ for definiteness)
preserves the boundary components of $\A$.  
If $u,v$ are lifts of $G$ and $H^2$ respectively to $\ti \A$ then 
$[u,v]$ has mean rotation zero and is a lift of $F^2 = [G,H^2].$
Again applying Theorem 2.1 of \cite{Fr6} we conclude that $F^2$ has
a fixed point in the interior of $U$.  

We want now to show this is also true in the case that both $G$ and
$H$ switch the boundary components of $\A$.  In that case $GH$ and $HG$
must preserve the boundary components of $\A$.  Let $g$ and $h$ be
lifts of $G$ and $H$ respectively to $\ti \A.$ They switch the ends of
$\ti \A$ and hence do not have mean rotation numbers.  But all elements
of the subgroup of the group generated by $g$ and $h$ consisting of
elements which preserve the ends of of $\ti \A$ will have
well defined mean rotation numbers.  This subgroup consists of all
elements which can be expressed as words of even length in $g$ and $h$.
Let $\rho_\mu(\alpha)$ denote the mean
rotation number of an element $\alpha$ in this subgroup and recall
$\rho_\mu$ is a homomorphism.  Then
\begin{align*}
\rho_\mu([g^2,h]) &= \rho_\mu(g^{-2} h^{-1} g^2 h) 
= \rho_\mu((g^{-2} h^{-1}g) gh) \\
&=
\rho_\mu(g^{-2} h^{-1}g^{-1} g^2) + \rho_\mu( gh)
= \rho_\mu( h^{-1}g^{-1} ) + \rho_\mu( gh) = 0.
\end{align*}
It follows from Theorem 2.1 of \cite{Fr6} again 
that $[g^2,h]$ has a fixed point, but it is a 
lift of $F^2 = [G^2,H]$, which must also have a fixed point.

Hence in all cases $\Fix(F^2)$ contains at least three points
and we can conclude from the previous case that
if $S \setminus \Fix(F^2)$ is non-empty, it
has negative Euler characteristic.
\end{proof}

\begin{lemma} \label{lem: fix=per}
Suppose $F\co S \to S$ is a homeomorphism,
  that each component $M$ of $S \setminus \Fix(F)$ has negative Euler characteristic and that
  for every $n>0$,    $F^n$  is isotopic to the identity rel $\Fix(F^n)$.  
Then $\Per(F) = \Fix(F).$
\end{lemma}

\begin{proof} 
Let  $f$ be the
restriction of $F$ to $M$.  We must show $\Per(f) = \emptyset.$  
Suppose to the contrary that $x \in \Per(f)$. Say it has period $p >1.$  Choose an arc $\alpha$   that initiates at $x$, terminates at  a point $y \in \Fix(F)$ and is otherwise disjoint from $\Fix(F)$.  By hypothesis, $F^p(\alpha)$ is isotopic to $\alpha$ relative to $\Fix(F^p)$ and hence relative to $\Fix(F) \cup \{x\}$.
Let $\ti f\co  \ti M \to \ti M$ be the identity lift and let $\ti \alpha_0 \subset H$ be a lift of $\alpha_0 = \alpha \setminus \{y\}$.  The initial endpoint $\ti x$ of $\ti \alpha$ is a lift of $x$ and the terminal end of $\ti \alpha_0$ converges to a point in $\sinfty$.  The isotopy of $F^p(\alpha)$ to $\alpha$ relative to $\Fix(F) \cup \{x\}$
lifts to an isotopy of $\ti f^p(\ti \alpha_0)$ to $\alpha_0$ relative to $\sinfty \cup \{\ti x\}$.  But this implies that $\ti f^p(\ti x) = \ti x$ in contradiction to the   
Brouwer plane translation theorem and that fact that $\ti f$ is fixed point free.   This contradicts the assumption 
that $\Per(f) \ne \emptyset.$
\end{proof}

\begin{prop}\label{prop: F = id}
Suppose $G,H \co S \to S$ are area preserving diffeomorphisms that are isotopic
to the identity and that commute with  
$F = [G,H]$.  Then $F^2 = id$.  
If  each component $M$ of $S \setminus \Fix(F)$ has negative Euler characteristic 
then $F = id$.
\end{prop}

\begin{proof}
We consider the second part first, so we assume $M$ has negative Euler characteristic.  Then according to Theorem 1.1 of \cite{fh:periodic}
either $F$ has points of arbitrarily high period or is the identity.
But according to Corollary~\ref{cor: commutator power}
 and
Lemma~\ref{lem: fix=per}, $F$ has no
points of period greater than one.  Hence it is the identity.

For the more general case we need only show that 
each component of $S \setminus \Fix(F^2)$ has negative Euler characteristic.
But this follows from Lemma~\ref{lem: hyperbolic complements}.
\end{proof}

\begin{ex}
Let $S^2$ be the unit sphere in $\R^3$ and let $G\co  S^2 \to S^2$ be rotation
through the angle $\pi$ around the $x$-axis.  Let 
$H\co  S^2 \to S^2$ be rotation
through the angle $\pi$ around an axis in the $xy$-plane which makes
an angle of $\pi/4$ with the $x$-axis.  Both $G$ and $H$ have 
order $2$.  One checks easily that
$F = [G,H]\co  S^2 \to S^2$ is rotation around the $z$-axis through
an angle of $\pi.$  Rotations through angle $\pi$ around perpendicular
axes commute.  Hence $F$ commutes with $G$ and $H$ and $F^2=id$
but $F \ne id$ which shows we cannot remove the exceptional case in
the preceding proposition.  The group generated by $G$ and $H$ is the dihedral 
group.
\end{ex}

We are now prepared to prove our main result.

\noindent
{\bf Theorem~\ref{no actions}}\qua {\sl Suppose $\G$ is an almost simple
group containing a subgroup isomorphic to the three-dimensional
integer Heisenberg group.  If $S$ is a closed oriented surface then every
homomorphism $\phi\co  \G \to \Diff_\omega(S)$ factors through a finite
group.}

\begin{proof}  As shown in the introduction we may assume that the action is by diffeomorphisms isotopic to the identity.
Since $\G$ contains a subgroup isomorphic to the three-dimensional 
integer Heisenberg group there are elements $G$ and $H$ such that $F = [G,H]$ 
has infinite order in $\G$ and commutes with $G$ and $H$.
By Proposition~\ref{prop: F = id} $\phi(F^2)$ is the identity.
It follows that $K = ker(\phi)$ has infinite order.  Since $\G$
is almost simple $K \subset \G$ has finite index.  Hence
$\phi$ factors through the finite group $\G/K.$
\end{proof}

\section{Nilpotent groups}

Suppose $\G$ is a finitely generated group and inductively define
$\G_i$ for $i > 0$ by $\G_0 = \G,\ \G_i = [\G, \G_{i-1}] :=$ the group
generated by $\{[g,h]\ |\  g \in \G,\ h \in \G_{i-1}\}.$  The group
$\G$ is called {\em nilpotent} if for some $n$,\ $\G_n = \{e\}.$
If $n$ is the smallest integer such that
$\G_n = \{e\}$ then $\G$ said to be of {\em nilpotent length}
$n$.  $\G$ is Abelian, if and only if its nilpotent length is 1.

\begin{thm} \label{nil}
Suppose that $\G$ is a finitely generated nilpotent subgroup 
of $\Diff_\omega(S)_0$.  If $S \ne S^2$ then $\G$ is Abelian;
if $S = S^2$ then $\G$ has an Abelian subgroup of index two.
\end{thm}

\begin{proof}
Suppose first that $S \ne S^2$.  
If $\G$ has nilpotent length $n > 1$ then $\G_{n-1}$
is non-trivial and its elements commute with all elements of $\G.$
Since $\G_{n-1}$ is generated by commutators 
there is a non-trivial element $F \in \G_{n-1}$ and
elements $G,H \in \G$ such that $F = [G,H].$  But Lemma~\ref{lem: hyperbolic complements} and 
Proposition~\ref{prop: F = id} then assert that $F = id$ contradicting
the assumption that $n > 1.$

In case $S = S^2$ we again assume 
$\G$ has nilpotent length $n > 1$ and hence that
there is a non-trivial element $F = [G,H] \in \G$ which commutes
with all elements of $\G.$
If $Fix(F)$ has at least three elements then again by Lemma~\ref{lem: hyperbolic complements} and Proposition~\ref{prop: F = id} 
$F = id$ contradicting the assumption that $n > 1$.
Hence we are reduced to the case $\Fix(F)$ consists of exactly
two points, $\{p,q\}.$

Let $\G'$ be subgroup of $\G$ consisting of all elements which
fix both $p$ and $q$.  Since every element of $\G$ commutes with
$F$ every element either fixes the two points $p$ and $q$ or
switches them.  Hence $\G'$ has index two in $\G$ and is nilpotent.

We will prove that $\G'$ is Abelian.  If not then then there is a
non-trivial element $F_0$ commuting will all of $\G'$ and elements
$G_0,H_0 \in \G'$ such that $F_0 = [G_0,H_0].$  Blow up the
two points $p$ and $q$ for the three diffeomorphisms
$F_0,G_0$ and $H_0$.  Since $F_0 = [G_0,H_0]$ the mean rotation
number of $F_0$ is zero.  It follows from Theorem~2.1 of \cite{Fr6}
that $F_0$ has an interior fixed point.  It then follows by the
argument above that $F_0 = id$ which is a contradiction.
\end{proof}

\begin{ex}  
One cannot replace $\Diff_\omega(S)_0$ with $\Diff_\omega(S)$ in the
preceding theorem. For example, there is an action of the
three-dimensional integer Heisenberg group on the two dimensional torus $T^2$
by area preserving diffeomorphisms. We will define the action first on
the universal cover $\mathbb R^2$. Choose an irrational number
$\alpha$ and define $\tilde G, \tilde H $ by $(x,y) \mapsto (x
+\alpha,y)$ and $(x,y) \mapsto (x,x+y)$ respectively.  The commutator
$\tilde F$ of $\tilde G$ and $\tilde H$ satisfies $(x,y) \mapsto (x,
y+\alpha)$ and so commutes with both $\tilde G$ and $\tilde H$.  The
maps $\tilde G$ and $\tilde H$ descend to area preserving
diffeomorphisms $G,H \co T^2 \to T^2$ that commute with their commutator.
It is easy to check that $G$ and $H$ generate a subgroup of
$\Diff_\omega(S)$ that is isomorphic to the three-dimensional
integer Heisenberg group. The map $H$ is a Dehn twist and so is not contained
in $\Diff_\omega(S)_0$.
\end{ex} 

A group $\G$ is {\it metabelian} if there is a homomorphism to an abelian group whose kernel is abelian.

\begin{cor} \label{metabelian} Any   finitely generated nilpotent subgroup $\G$
of $\Diff_\omega(S)$ has a finite index metabelian subgroup.  
\end{cor}

\begin{proof} 
If $S= S^2$, then $\G$ has an index two subgroup that is contained in
$\Diff_\omega(S)_0$ so Theorem~\ref{nil} implies that $\G$ is
virtually abelian.  For $S \ne S^2$ let $\Phi \co  \Diff_\omega(S) \to
MCG(S)$ be the natural map. By \cite{BLM}, $\Phi(\G)$ is virtually
abelian.  Thus $\G$ has a finite index subgroup $\G_0$ whose image
$\Phi(\G_0)$ is abelian.  Theorem~\ref{nil} implies that the kernel of
$\Phi|_{G_0}$ is abelian.  
\end{proof}

\rk{Acknowledgements}

John Franks is supported in part by NSF grant 
DMS0099640.\nl Michael Handel is supported in part by NSF grant DMS0103435.

\end{document}